# Two Series Formulae from the Solution of a System of Diffusion Equations

SIDDHARTH. G. CHATTERJEE

Department of Paper and Bioprocess Engineering
SUNY College of Environmental Science and Forestry
1 Forestry Drive, Syracuse, New York 13210, USA
E-mail: schatterjee@esf.edu

**Abstract.** From the solution of the heat and mass diffusion equations that describe the exothermic physical absorption of a gas into a liquid, compact formulae for the sums of two infinite series are derived. The method is general and should be applicable to other systems of linear, parabolic partial differential equations.

## 1. Introduction

The main objective of this brief paper is to present a method of deriving compact, closed-form formulae for the sums of infinite series from the solution of a system of linear, parabolic partial differential equations.

## 2. Theory

Consider the nonisothermal physical absorption of a gas into a liquid. The governing (dimensionless) equations of the film-penetration model that describe this phenomenon are [1]:

$$\frac{\partial a}{\partial \tau} = \frac{\partial^2 a}{\partial x^2} \tag{1}$$

$$\varepsilon \frac{\partial \theta}{\partial \tau} = \frac{\partial^2 \theta}{\partial x^2} \tag{2}$$

with

$$a = 0 \text{ and } \theta = 0 \text{ at } \tau = 0 \text{ for } 0 \leq x \leq 1 \tag{3}$$

$$a = 0 \text{ and } \theta = 0 \text{ for } \tau > 0 \text{ at } x = 1 \tag{4}$$

$$a = w_1 \theta + w_2 \text{ for } \tau > 0 \text{ at } x = 0 \tag{5}$$

$$\frac{\partial \theta}{\partial x} = \beta_s \left(\frac{\partial a}{\partial x}\right) \text{ for } \tau > 0 \text{ at } x = 0 \tag{6}$$

In the above, $\theta$ and $a$ are the dimensionless temperature of the liquid and concentration of the dissolved gas at (dimensionless) location $x$ and (dimensionless) time $\tau$, whereas $w_1$, $w_2$, $\beta_s$, and $\varepsilon$ are dimensionless constants whose definitions can be found in [1].

(1) and (2), subject to (3) – (6), can be solved by taking Laplace transforms (with respect to $\tau$) and inverting the resulting transforms by the theorem of residues [2]. The Laplace transforms of $a$ and $\theta$ are given by [1]:

$$\bar{a}(x,p) = \frac{w_2[\cosh(x\sqrt{p}) - \coth\sqrt{p}\,\sinh(x\sqrt{p})]}{p\left[1 - \left(\frac{w_1\beta_s}{\sqrt{\varepsilon}}\right)\coth\sqrt{p}\,\tanh(\sqrt{\varepsilon p})\right]} \qquad (7)$$

$$\bar{\theta}(x,p) = \frac{w_2\beta_s[\cosh(x\sqrt{\varepsilon p}) - \coth(\sqrt{\varepsilon p})\sinh(x\sqrt{\varepsilon p})]}{p\sqrt{\varepsilon}\left[\coth(\sqrt{\varepsilon p})\tanh\sqrt{p} - \left(\frac{w_1\beta_s}{\sqrt{\varepsilon}}\right)\right]} \qquad (8)$$

where $p$ is the Laplace transform parameter, which has the physical significance of being a dimensionless surface-renewal rate (a positive quantity) in the film-penetration model. The inversion of (7) and (8) yields [1]

$$a(x,\tau) = \frac{w_2(1-x)}{1 - w_1\beta_s} + 2w_2\sum_{n=1}^{\infty}\frac{\cos(q_n\sqrt{\varepsilon})\sin[q_n(1-x)]\exp(-q_n^2\tau)}{q_n G(q_n)} \qquad (9)$$

$$\theta(x,\tau) = \frac{w_2\beta_s(1-x)}{1 - w_1\beta_s} + \frac{2w_2\beta_s}{\sqrt{\varepsilon}}\sum_{n=1}^{\infty}\frac{\cos q_n\,\sin[q_n\sqrt{\varepsilon}(1-x)]\exp(-q_n^2\tau)}{q_n G(q_n)} \qquad (10)$$

where

$$G(q_n) = (1 - w_1\beta_s)\cos(q_n\sqrt{\varepsilon})\cos q_n - \left(\sqrt{\varepsilon} - \frac{w_1\beta_s}{\sqrt{\varepsilon}}\right)\sin(q_n\sqrt{\varepsilon})\sin q_n \qquad (11)$$

and the $q_n$'s are the positive roots of the equation

$$\tan q_n = \frac{w_1\beta_s}{\sqrt{\varepsilon}}\tan(q_n\sqrt{\varepsilon}); \; n = 1, 2, 3, \ldots, \qquad (12)$$

Taking the Laplace transform of the surface temperature $\theta(x = 0, \tau)$ obtained from (10) and equating the resulting expression to $\bar{\theta}(x = 0, p)$ got from (8) yields the first formula:

$$\frac{\sqrt{\varepsilon}}{1-\gamma} + 2\sum_{n=1}^{\infty}\frac{\cos q_n\,\sin(q_n\sqrt{\varepsilon})}{q_n(1 + q_n^2/p)G(q_n)} = \frac{1}{\coth(\sqrt{\varepsilon p})\tanh\sqrt{p} - \frac{\gamma}{\sqrt{\varepsilon}}} \qquad (13)$$

where

$$\gamma = w_1\beta_s \qquad (14)$$

Evaluating $(\partial a/\partial x)_{x=0}$ from (9), taking its Laplace transform, and equating the resulting expression to $(\partial \bar{a}/\partial x)_{x=0}$ obtained from (7) gives the second formula:

$$\frac{1}{1-\gamma} + 2\sum_{n=1}^{\infty}\frac{\cos(q_n\sqrt{\varepsilon})\cos q_n}{(1 + q_n^2/p)G(q_n)} = \frac{\sqrt{p}\,\coth\sqrt{p}}{1 - \frac{\gamma}{\sqrt{\varepsilon}}\coth\sqrt{p}\,\tanh(\sqrt{\varepsilon p})} \qquad (15)$$

The numerator in the fraction on the right-hand-side (RHS) of (15) can be expressed as the sum of an infinite series that can be obtained from the solution of the film-penetration model for isothermal absorption of a gas into a liquid [3]. The governing equations of this model in dimensionless form may be written as:

$$\frac{\partial a}{\partial \tau} = \frac{\partial^2 a}{\partial x^2} \qquad (16)$$

with

$$a = 0 \text{ at } \tau = 0 \text{ for } 0 \leq x \leq 1 \tag{17}$$

$$a = 0 \text{ for } \tau > 0 \text{ at } x = 1 \tag{18}$$

$$a = w_2 \text{ for } \tau > 0 \text{ at } x = 0 \tag{19}$$

By taking the Laplace transform of (16) and using (17) – (19), it can be shown that

$$\bar{a}(x,p) = \frac{w_2[\cosh(x\sqrt{p}) - \coth\sqrt{p} \sinh(x\sqrt{p})]}{p} \tag{20}$$

which expression can also be obtained by setting $\gamma = 0$ (i.e., $\beta_s = 0$) in (7). Equating the Laplace transform of $(\partial a/\partial x)_{x=0}$ obtained from [3] to $(\partial \bar{a}/\partial x)_{x=0}$ obtained from (20) yields the third formula:

$$1 + 2 \sum_{n=1}^{\infty} \frac{1}{\left(1 + \frac{n^2 \pi^2}{p}\right)} = \sqrt{p} \coth\sqrt{p} \tag{21}$$

which is a known identity [4, 5]. Setting $\gamma = 0$ in (12), (11) and (15) gives the following:

$$\tan q_n = 0, \text{ i. e.}, q_n = n\pi; \quad n = 1, 2, 3, ...., \tag{22}$$

$$G(q_n) = \cos(q_n \sqrt{\varepsilon}) \cos q_n \tag{23}$$

$$1 + 2 \sum_{n=1}^{\infty} \frac{\cos(q_n \sqrt{\varepsilon}) \cos q_n}{(1 + q_n^2/p) G(q_n)} = \sqrt{p} \coth\sqrt{p} \tag{24}$$

Substituting (22) and (23) into (24) yields (21), which is, therefore, a special case of (15). Setting $p = \pi^2$ in (21) yields the following interesting result:

$$1 + 2 \sum_{n=1}^{\infty} \frac{1}{1 + n^2} = \pi \coth \pi \tag{25}$$

It is to be noted that the series on the left-hand-side (LHS) of (21), and also those in (13) and (15), are convergent for all values of $p < \infty$ (assuming $p$ to be real). A formula for the hyperbolic cotangent is available [6], which can be written as

$$1 + \frac{p}{3} - \frac{p^2}{45} + \frac{2}{945} p^3 - \cdots + \frac{2^{2n} B_{2n}}{(2n)!} p^n - \cdots = \sqrt{p} \coth\sqrt{p} \tag{26}$$

(26) is valid for $p < \pi^2$ with $B_n$ being the $n$th Bernoulli number.

## 3. Results and Discussion

(11) – (15) reveal that the only two fundamental dimensionless constants are $\gamma$ and $\varepsilon$. Table 1 reports the first twenty five roots of (12) for $\gamma = -0.03421$ and $\varepsilon = 2.64489 \times 10^{-3}$, which represent the water vapor–lithium bromide system [1]. For values of $p$ ranging over ten orders of magnitude, Tables 2 and 3 compare the left-hand and right-hand sides of (13) and (15), respectively, for the above values of $\gamma$ and $\varepsilon$ using the roots in Table 1. The maximum discrepancy between the LHS and RHS of (13) in Table 2 is about 2% for $p = 10^5 - 10^6$. In contrast, Table 3 shows that the difference between the LHS and RHS of (15) increases rapidly for $p > 100$. This is due to a finite number of terms (i.e., 25) used to represent the infinite series in (15) as was discovered when the number of terms was reduced to 20, which further widened the discrepancy between the LHS and RHS of (15).

We conclude by noting that for single, linear parabolic or elliptic partial differential equations, formulae for the sums of infinite series can also be obtained by the method of finite integral transforms (taken with respect to the spatial coordinate) as indicated by Özişik [7] in the solution of Laplace's equation over a finite rectangle from which the following formula can be derived:

$$\frac{2}{b}\sum_{n=1}^{\infty}\frac{\vartheta_n}{\sigma_m^2+\vartheta_n^2}\sin(\vartheta_n y) = \frac{\sinh[\sigma_m(b-y)]}{\sinh(\sigma_m b)} \qquad (27)$$

$$\sigma_m = \frac{m\pi}{a}, \qquad \vartheta_n = \frac{n\pi}{b}, \qquad m,n = 1,2,3,\dots., \qquad (28)$$

where *a*, *b* and *y* are real. However, this method hinges upon the availability of an appropriate kernel, which will depend upon the nature of the boundary conditions of the partial differential equation.

**Acknowledgment**

The author is grateful to Mr. Ross David Mazur and Dr. Wayne S. Amato for pointing out [4] and (5), respectively.

TABLE 1

First 25 Roots of (12) for $\gamma = -0.03421$ and $\varepsilon = 2.64489 \times 10^{-3}$.[a]

| n | $q_n$ | residual (%)[b] |
|---|---|---|
| 1 | 3.0371 | 0.0993 |
| 2 | 6.0715 | 0.0991 |
| 3 | 9.1001 | 0.0999 |
| 4 | 12.1199 | 0.0951 |
| 5 | 15.1273 | 0.1032 |
| 6 | 18.1187 | 0.1182 |
| 7 | 21.0910 | 0.1262 |
| 8 | 24.0425 | 0.1620 |
| 9 | 26.9750 | 0.2372 |
| 10 | 29.8947 | 1.0935 |
| 11 | 32.8119 | 0.2518 |
| 12 | 35.7377 | 0.0904 |
| 13 | 38.6803 | 0.0664 |
| 14 | 41.6435 | 0.0461 |
| 15 | 44.6269 | 0.0287 |
| 16 | 47.6277 | 0.0399 |
| 17 | 50.6424 | 0.0836 |
| 18 | 53.6676 | 0.1187 |
| 19 | 56.6998 | 0.2546 |
| 20 | 59.7361 | 0.8836 |
| 21 | 62.7735 | 0.9581 |
| 22 | 65.8094 | 0.4222 |
| 23 | 68.8409 | 0.3212 |
| 24 | 71.8650 | 0.2747 |
| 25 | 74.8783 | 0.2597 |

[a] The first 20 roots were reported in [1].
[b] Defined as $\left|[\tan q_n - (\gamma/\sqrt{\varepsilon})\tan(q_n\sqrt{\varepsilon})]/\tan q_n\right| \times 100$.

TABLE 2
*Comparison of the LHS and RHS of formula no. 1 (13) for $\gamma = -0.03421$, $\varepsilon = 2.64489 \times 10^{-3}$ using the roots in Table 1.*

| $p$ | LHS [equation (13)] | RHS [equation (13)] |
|---|---|---|
| $10^{-4}$ | 0.04973 | 0.04973 |
| $10^{-3}$ | 0.04974 | 0.04974 |
| $10^{-2}$ | 0.04989 | 0.04989 |
| $10^{-1}$ | 0.05131 | 0.05131 |
| $10^{0}$ | 0.06457 | 0.06457 |
| $10^{1}$ | 0.14602 | 0.14607 |
| $10^{2}$ | 0.35950 | 0.35995 |
| $10^{3}$ | 0.56938 | 0.57284 |
| $10^{4}$ | 0.58911 | 0.60049 |
| $10^{5}$ | 0.58757 | 0.60051 |
| $10^{6}$ | 0.58761 | 0.60051 |

TABLE 3
*Comparison of the LHS and RHS of formula No. 2 (15) for $\gamma = -0.03421$, $\varepsilon = 2.64489 \times 10^{-3}$ using the roots in Table 1.*

| $p$ | LHS [equation (15)] | RHS [equation (15)] |
|---|---|---|
| $10^{-4}$ | 0.96695 | 0.96695 |
| $10^{-3}$ | 0.96723 | 0.96723 |
| $10^{-2}$ | 0.96999 | 0.97003 |
| $10^{-1}$ | 0.99739 | 0.99785 |
| $10^{0}$ | 1.25209 | 1.25663 |
| $10^{1}$ | 2.81978 | 2.86524 |
| $10^{2}$ | 7.15260 | 7.60546 |
| $10^{3}$ | 15.22671 | 19.57211 |
| $10^{4}$ | 26.77400 | 60.05292 |
| $10^{5}$ | 31.43094 | 189.89883 |
| $10^{6}$ | 32.06626 | 600.51281 |